# CURRENT STATUS DATA WITH COMPETING RISKS: LIMITING DISTRIBUTION OF THE MLE

By Piet Groeneboom, Marloes H. Maathuis[1]
and Jon A. Wellner[2]

*Delft University of Technology and Vrije Universiteit Amsterdam,
University of Washington and University of Washington*

We study nonparametric estimation for current status data with competing risks. Our main interest is in the nonparametric maximum likelihood estimator (MLE), and for comparison we also consider a simpler "naive estimator." Groeneboom, Maathuis and Wellner [*Ann. Statist.* (2008) **36** 1031–1063] proved that both types of estimators converge globally and locally at rate $n^{1/3}$. We use these results to derive the local limiting distributions of the estimators. The limiting distribution of the naive estimator is given by the slopes of the convex minorants of correlated Brownian motion processes with parabolic drifts. The limiting distribution of the MLE involves a new self-induced limiting process. Finally, we present a simulation study showing that the MLE is superior to the naive estimator in terms of mean squared error, both for small sample sizes and asymptotically.

**1. Introduction.** We study nonparametric estimation for current status data with competing risks. The set-up is as follows. We analyze a system that can fail from $K$ competing risks, where $K \in \mathbb{N}$ is fixed. The random variables of interest are $(X, Y)$, where $X \in \mathbb{R}$ is the failure time of the system, and $Y \in \{1, \ldots, K\}$ is the corresponding failure cause. We cannot observe $(X, Y)$ directly. Rather, we observe the "current status" of the system at a single random observation time $T \in \mathbb{R}$, where $T$ is independent of $(X, Y)$. This means that at time $T$, we observe whether or not failure occurred, and if and only if failure occurred, we also observe the failure cause $Y$. Such

Received September 2006; revised April 2007.
[1]Supported in part by NSF Grant DMS-02-03320.
[2]Supported in part by NSF Grants DMS-02-03320 and DMS-05-03822 and by NI-AID Grant 2R01 AI291968-04.

data arise naturally in cross-sectional studies with several failure causes, and generalizations arise in HIV vaccine clinical trials (see [10]).

We study nonparametric estimation of the sub-distribution functions $F_{01}, \ldots, F_{0K}$, where $F_{0k}(s) = P(X \leq s, Y = k)$, $k = 1, \ldots, K$. Various estimators for this purpose were introduced in [10, 12], including the nonparametric maximum likelihood estimator (MLE), which is our primary focus. For comparison we also consider the "naive estimator," an alternative to the MLE discussed in [12]. Characterizations, consistency and $n^{1/3}$ rates of convergence of these estimators were established in Groeneboom, Maathuis and Wellner [8]. In the current paper we use these results to derive the local limiting distributions of the estimators.

1.1. *Notation.* The following notation is used throughout. The observed data are denoted by $(T, \Delta)$, where $T$ is the observation time and $\Delta = (\Delta_1, \ldots, \Delta_{K+1})$ is an indicator vector defined by $\Delta_k = 1\{X \leq T, Y = k\}$ for $k = 1, \ldots, K$, and $\Delta_{K+1} = 1\{X > T\}$. Let $(T_i, \Delta^i)$, $i = 1, \ldots, n$, be $n$ i.i.d. observations of $(T, \Delta)$, where $\Delta^i = (\Delta_1^i, \ldots, \Delta_{K+1}^i)$. Note that we use the superscript $i$ as the index of an observation, and not as a power. The order statistics of $T_1, \ldots, T_n$ are denoted by $T_{(1)}, \ldots, T_{(n)}$. Furthermore, $G$ is the distribution of $T$, $G_n$ is the empirical distribution of $T_i$, $i = 1, \ldots, n$, and $\mathbb{P}_n$ is the empirical distribution $(T_i, \Delta^i)$, $i = 1, \ldots, n$. For any vector $(x_1, \ldots, x_K) \in \mathbb{R}^K$ we define $x_+ = \sum_{k=1}^K x_k$, so that, for example, $\Delta_+ = \sum_{k=1}^K \Delta_k$ and $F_{0+}(s) = \sum_{k=1}^K F_{0k}(s)$. For any $K$-tuple $F = (F_1, \ldots, F_K)$ of sub-distribution functions, we define $F_{K+1}(s) = \int_{u>s} dF_+(u) = F_+(\infty) - F_+(s)$.

We denote the right-continuous derivative of a function $f : \mathbb{R} \mapsto \mathbb{R}$ by $f'$ (if it exists). For any function $f : \mathbb{R} \mapsto \mathbb{R}$, we define the *convex minorant* of $f$ to be the largest convex function that is pointwise bounded by $f$. For any interval $I$, $D(I)$ denotes the collection of cadlag functions on $I$. Finally, we use the following definition for integrals and indicator functions:

DEFINITION 1.1. Let $dA$ be a Lebesgue–Stieltjes measure, and let $W$ be a Brownian motion process. For $t < t_0$, we define $1_{[t_0,t)}(u) = -1_{[t,t_0)}(u)$,

$$\int_{[t_0,t)} f(u) \, dA(u) = -\int_{[t,t_0)} f(u) \, dA(u)$$

and

$$\int_{t_0}^t f(u) \, dW(u) = -\int_t^{t_0} f(u) \, dW(u).$$

1.2. *Assumptions.* We prove the local limiting distributions of the estimators at a fixed point $t_0$, under the following assumptions: (a) The observation time $T$ is independent of the variables of interest $(X, Y)$. (b) For each $k = 1, \ldots, K$, $0 < F_{0k}(t_0) < F_{0k}(\infty)$, and $F_{0k}$ and $G$ are continuously differentiable at $t_0$ with positive derivatives $f_{0k}(t_0)$ and $g(t_0)$. (c) The system cannot fail from two or more causes at the same time. Assumptions (a) and (b) are essential for the development of the theory. Assumption (c) ensures that the failure cause is well defined. This assumption is always satisfied by defining simultaneous failure from several causes as a new failure cause.

1.3. *The estimators.* We first consider the MLE. The MLE $\widehat{F}_n = (\widehat{F}_{n1}, \ldots, \widehat{F}_{nK})$ is defined by $l_n(\widehat{F}_n) = \max_{F \in \mathcal{F}_K} l_n(F)$, where

$$(1) \quad l_n(F) = \int \left\{ \sum_{k=1}^{K} \delta_k \log F_k(t) + (1 - \delta_+) \log(1 - F_+(t)) \right\} d\mathbb{P}_n(t, \delta),$$

and $\mathcal{F}_K$ is the collection of $K$-tuples $F = (F_1, \ldots, F_K)$ of sub-distribution functions on $\mathbb{R}$ with $F_+ \leq 1$. The naive estimator $\widetilde{F}_n = (\widetilde{F}_{n1}, \ldots, \widetilde{F}_{nK})$ is defined by $l_{nk}(\widetilde{F}_{nk}) = \max_{F_k \in \mathcal{F}} l_{nk}(F_k)$, for $k = 1, \ldots, K$, where $\mathcal{F}$ is the collection of distribution functions on $\mathbb{R}$, and

$$(2) \quad l_{nk}(F_k) = \int \{\delta_k \log F_k(t) + (1 - \delta_k) \log(1 - F_k(t))\} d\mathbb{P}_n(t, \delta),$$

$$k = 1, \ldots, K.$$

Note that $\widetilde{F}_{nk}$ only uses the $k$th entry of the $\Delta$-vector, and is simply the MLE for the reduced current status data $(T, \Delta_k)$. Thus, the naive estimator splits the optimization problem into $K$ separate well-known problems. The MLE, on the other hand, estimates $F_{01}, \ldots, F_{0K}$ simultaneously, accounting for the fact that $\sum_{k=1}^{K} F_{0k}(s) = P(X \leq s)$ is the overall failure time distribution. This relation is incorporated both in the object function $l_n(F)$ [via the term $\log(1 - F_+)$] and in the space $\mathcal{F}_K$ over which $l_n(F)$ is maximized (via the constraint $F_+ \leq 1$).

1.4. *Main results.* The main results in this paper are the local limiting distributions of the MLE and the naive estimator. The limiting distribution of $\widetilde{F}_{nk}$ corresponds to the limiting distribution of the MLE for the reduced current status data $(T, \Delta_k)$. Thus, it is given by the slope of the convex minorant of a two-sided Brownian motion process plus parabolic drift ([9], Theorem 5.1, page 89) known as Chernoff's distribution. The joint limiting distribution of $(\widetilde{F}_{n1}, \ldots, \widetilde{F}_{nK})$ follows by noting that the $K$ Brownian motion processes have a multinomial covariance structure, since $\Delta | T \sim \text{Mult}_{K+1}(1, (F_{01}(T), \ldots, F_{0,K+1}(T)))$. The drifted Brownian motion processes and their convex minorants are specified in Definitions 1.2 and 1.5. The limiting distribution of the naive estimator is given in Theorem 1.6, and is simply a $K$-dimensional version of the limiting distribution for current status data. A formal proof of this result can be found in [14], Section 6.1.

DEFINITION 1.2. Let $W = (W_1, \ldots, W_K)$ be a $K$-tuple of two-sided Brownian motion processes originating from zero, with mean zero and covariances

$$(3) \quad E\{W_j(t)W_k(s)\} = (|s| \wedge |t|)1\{st > 0\}\Sigma_{jk}, \qquad s, t \in \mathbb{R}, j, k \in \{1, \ldots, K\},$$

where $\Sigma_{jk} = g(t_0)^{-1}\{1\{j = k\}F_{0k}(t_0) - F_{0j}(t_0)F_{0k}(t_0)\}$. Moreover, $V = (V_1, \ldots, V_K)$ is a vector of drifted Brownian motions, defined by

$$(4) \quad V_k(t) = W_k(t) + \tfrac{1}{2}f_{0k}(t_0)t^2, \qquad k = 1, \ldots, K.$$

Following the convention introduced in Section 1.1, we write $W_+ = \sum_{k=1}^{K} W_k$ and $V_+ = \sum_{k=1}^{K} V_k$. Finally, we use the shorthand notation $a_k = (F_{0k}(t_0))^{-1}$, $k = 1, \ldots, K+1$.

REMARK 1.3. Note that $W$ is the limit of a rescaled version of $W_n = (W_{n1}, \ldots, W_{nK})$, and that $V$ is the limit of a recentered and rescaled version of $V_n = (V_{n1}, \ldots, V_{nK})$, where $W_{nk}$ and $V_{nk}$ are defined by (17) and (6) of [8]:

$$(5) \quad \begin{aligned} W_{nk}(t) &= \int_{u \leq t} \{\delta_k - F_{0k}(t_0)\} d\mathbb{P}_n(u, \delta), & t \in \mathbb{R}, k = 1, \ldots, K, \\ V_{nk}(t) &= \int_{u \leq t} \delta_k \, d\mathbb{P}_n(u, \delta), & t \in \mathbb{R}, k = 1, \ldots, K. \end{aligned}$$

REMARK 1.4. We define the correlation between Brownian motions $W_j$ and $W_k$ by

$$r_{jk} = \frac{\Sigma_{jk}}{\sqrt{\Sigma_{jj}\Sigma_{kk}}} = -\frac{\sqrt{F_{0j}(t_0)F_{0k}(t_0)}}{\sqrt{(1 - F_{0j}(t_0))(1 - F_{0k}(t_0))}}.$$

Thus, the Brownian motions are negatively correlated, and this negative correlation becomes stronger as $t_0$ increases. In particular, it follows that $r_{12} \to -1$ as $F_{0+}(t_0) \to 1$, in the case of $K = 2$ competing risks.

DEFINITION 1.5. Let $\widetilde{H} = (\widetilde{H}_1, \ldots, \widetilde{H}_K)$ be the vector of convex minorants of $V$, that is, $\widetilde{H}_k$ is the convex minorant of $V_k$, for $k = 1, \ldots, K$. Let $\widetilde{F} = (\widetilde{F}_1, \ldots, \widetilde{F}_K)$ be the vector of right derivatives of $\widetilde{H}$.

THEOREM 1.6. *Under the assumptions of Section 1.2,*

$$n^{1/3}\{\widetilde{F}_n(t_0 + n^{-1/3}t) - F_0(t_0)\} \to_d \widetilde{F}(t)$$

*in the Skorohod topology on* $(D(\mathbb{R}))^K$.

The limiting distribution of the MLE is given by the slopes of a new self-induced process $\widehat{H} = (\widehat{H}_1, \ldots, \widehat{H}_K)$, defined in Theorem 1.7. We say that the process $\widehat{H}$ is "self-induced," since each component $\widehat{H}_k$ is defined in terms of the other components through $\widehat{H}_+ = \sum_{j=1}^K \widehat{H}_j$. Due to this self-induced nature, existence and uniqueness of $\widehat{H}$ need to be formally established (Theorem 1.7). The limiting distribution of the MLE is given in Theorem 1.8. These results are proved in the remainder of the paper.

THEOREM 1.7. *There exists an almost surely unique $K$-tuple $\widehat{H} = (\widehat{H}_1, \ldots, \widehat{H}_K)$ of convex functions with right-continuous derivatives $\widehat{F} = (\widehat{F}_1, \ldots, \widehat{F}_K)$, satisfying the following three conditions:*

(i) $a_k \widehat{H}_k(t) + a_{K+1} \widehat{H}_+(t) \leq a_k V_k(t) + a_{K+1} V_+(t)$, *for* $k = 1, \ldots, K$, $t \in \mathbb{R}$.

(ii) $\int \{a_k \widehat{H}_k(t) + a_{K+1} \widehat{H}_+(t) - a_k V_k(t) - a_{K+1} V_+(t)\} d\widehat{F}_k(t) = 0$, $k = 1, \ldots, K$.

(iii) *For all $M > 0$ and $k = 1, \ldots, K$, there are points $\tau_{1k} < -M$ and $\tau_{2k} > M$ so that $a_k \widehat{H}_k(t) + a_{K+1} \widehat{H}_+(t) = a_k V_k(t) + a_{K+1} V_+(t)$ for $t = \tau_{1k}$ and $t = \tau_{2k}$.*

THEOREM 1.8. *Under the assumptions of Section* 1.2,
$$n^{1/3}\{\widehat{F}_n(t_0 + n^{-1/3}t) - F_0(t_0)\} \to_d \widehat{F}(t)$$
*in the Skorohod topology on* $(D(\mathbb{R}))^K$.

Thus, the limiting distributions of the MLE and the naive estimator are given by the slopes of the limiting processes $\widehat{H}$ and $\widetilde{H}$, respectively. In order to compare $\widehat{H}$ and $\widetilde{H}$, we note that the convex minorant $\widetilde{H}_k$ of $V_k$ can be defined as the almost surely unique convex function $\widetilde{H}_k$ with right-continuous derivative $\widetilde{F}_k$ that satisfies: (i) $\widetilde{H}_k(t) \leq V_k(t)$ for all $t \in \mathbb{R}$, and (ii) $\int \{\widetilde{H}_k(t) - \widetilde{V}_k(t)\} d\widetilde{F}_k(t) = 0$. Comparing this to the definition of $\widehat{H}_k$ in Theorem 1.7, we see that the definition of $\widehat{H}_k$ contains the extra terms $\widehat{H}_+$ and $V_+$, which come from the term $\log(1 - F_+(t))$ in the log likelihood (1). The presence of $\widehat{H}_+$ in Theorem 1.7 causes $\widehat{H}$ to be self-induced. In contrast, the processes $\widetilde{H}_k$ for the naive estimator depend only on $V_k$, so that $\widetilde{H}$ is not self-induced. However, note that the processes $\widetilde{H}_1, \ldots, \widetilde{H}_K$ are correlated, since the Brownian motions $W_1, \ldots, W_K$ are correlated (see Definition 1.2).

1.5. *Outline.* This paper is organized as follows. In Section 2 we discuss the new self-induced limiting processes $\widehat{H}$ and $\widehat{F}$. We give various interpretations of these processes and prove the uniqueness part of Theorem 1.7. Section 3 establishes convergence of the MLE to its limiting distribution (Theorem 1.8). Moreover, in this proof we automatically obtain existence of $\widehat{H}$, hence completing the proof of Theorem 1.7. This approach to proving existence of the limiting processes is different from the one followed by [5, 6] for the estimation of convex functions, who establish existence and uniqueness of the limiting process before proving convergence. In Section 4 we compare the estimators in a simulation study, and show that the MLE is superior to the naive estimator in terms of mean squared error, both for small sample sizes and asymptotically. We also discuss computation of the estimators in Section 4. Technical proofs are collected in Section 5.

**2. Limiting processes.** We now discuss the new self-induced processes $\widehat{H}$ and $\widehat{F}$ in more detail. In Section 2.1 we give several interpretations of these processes, and illustrate them graphically. In Section 2.2 we prove tightness of $\{\widehat{F}_k - f_{0k}(t_0)t\}$ and $\{\widehat{H}_k(t) - V_k(t)\}$, for $t \in \mathbb{R}$. These results are used in Section 2.3 to prove almost sure uniqueness of $\widehat{H}$ and $\widehat{F}$.

2.1. *Interpretations of $\widehat{H}$ and $\widehat{F}$.* Let $k \in \{1, \ldots, K\}$. Theorem 1.7(i) and the convexity of $\widehat{H}_k$ imply that $a_k \widehat{H}_k + a_{K+1} \widehat{H}_+$ is a convex function that lies below $a_k V_k + a_{K+1} V_+$. Hence, $a_k \widehat{H}_k + a_{K+1} \widehat{H}_+$ is bounded above by the convex minorant of $a_k V_k + a_{K+1} V_+$. This observation leads directly to the following proposition about the points of touch between $a_k \widehat{H}_k + a_{K+1} \widehat{H}_+$ and $a_k V_k + a_{K+1} V_+$:

PROPOSITION 2.1. *For each $k = 1, \ldots, K$, we define $\mathcal{N}_k$ and $\widehat{\mathcal{N}}_k$ by*

(6) $\mathcal{N}_k = \{points\ of\ touch\ between\ a_k V_k + a_{K+1} V_+\ and\ its\ convex\ minorant\},$

(7) $\widehat{\mathcal{N}}_k = \{points\ of\ touch\ between\ a_k V_k + a_{K+1} V_+\ and\ a_k \widehat{H}_k + a_{K+1} \widehat{H}_+\}$

*Then the following properties hold:* (i) $\widehat{\mathcal{N}}_k \subset \mathcal{N}_k$, *and* (ii) *At points $t \in \widehat{\mathcal{N}}_k$, the right and left derivatives of $a_k \widehat{H}_k(t) + a_{K+1} \widehat{H}_+(t)$ are bounded above and below by the right and left derivatives of the convex minorant of $a_k V_k(t) + a_{K+1} V_+(t)$.*

Since $a_k V_k + a_{K+1} V_+$ is a Brownian motion process plus parabolic drift, the point process $\mathcal{N}_k$ is well known from [4]. On the other hand, little is known about $\widehat{\mathcal{N}}_k$, due to the self-induced nature of this process. However, Proposition 2.1(i) relates $\widehat{\mathcal{N}}_k$ to $\mathcal{N}_k$, and this allows us to deduce properties of $\widehat{\mathcal{N}}_k$ and the associated processes $\widehat{H}_k$ and $\widehat{F}_k$. In particular, Proposition 2.1(i) implies that $\widehat{F}_k$ is piecewise constant, and that $\widehat{H}_k$ is piecewise linear (Corollary 2.2). Moreover, Proposition 2.1(i) is essential for the proof of Proposition 2.16, where it is used to establish expression (30). Proposition 2.1(ii) is not used in the sequel.

COROLLARY 2.2. *For each $k \in \{1,\ldots,K\}$, the following properties hold almost surely:* (i) $\widehat{\mathcal{N}}_k$ *has no condensation points in a finite interval, and* (ii) $\widehat{F}_k$ *is piecewise constant and $\widehat{H}_k$ is piecewise linear.*

PROOF. $\mathcal{N}_k$ is a stationary point process which, with probability 1, has no condensation points in a finite interval (see [4]). Together with Proposition 2.1(i), this yields that with probability 1, $\widehat{\mathcal{N}}_k$ has no condensation points in a finite interval. Conditions (i) and (ii) of Theorem 1.7 imply that $\widehat{F}_k$ can only increase at points $t \in \mathcal{N}_k$. Hence, $\widehat{F}_k$ is piecewise constant and $\widehat{H}_k$ is piecewise linear. $\square$

Thus, conditions (i) and (ii) of Theorem 1.7 imply that $a_k \widehat{H}_k + a_{K+1} \widehat{H}_+$ is a piecewise linear convex function, lying below $a_k V_k + a_{K+1} V_+$, and touching $a_k V_k + a_{K+1} V_+$ whenever $\widehat{F}_k$ jumps. We illustrate these processes using the following example with $K = 2$ competing risks:

EXAMPLE 2.3. Let $K = 2$, and let $T$ be independent of $(X, Y)$. Let $T$, $Y$ and $X|Y$ be distributed as follows: $G(t) = 1 - \exp(-t)$, $P(Y = k) = k/3$ and $P(X \leq t | Y = k) = 1 - \exp(-kt)$ for $k = 1, 2$. This yields $F_{0k}(t) = (k/3)\{1 - \exp(-kt)\}$ for $k = 1, 2$.

Figure 1 shows the limiting processes $a_k V_k + a_{K+1} V_+$, $a_k \widehat{H}_k + a_{K+1} \widehat{H}_+$, and $\widehat{F}_k$, for this model with $t_0 = 1$. The relevant parameters at the point $t_0 = 1$ are

$$F_{01}(1) = 0.21, \qquad F_{02}(1) = 0.58,$$
$$f_{01}(1) = 0.12, \qquad f_{02}(1) = 0.18, \qquad g(1) = 0.37.$$

The processes shown in Figure 1 are approximations, obtained by computing the MLE for sample size $n = 100{,}000$ (using the algorithm described in Section 4), and then computing the localized processes $V_{nk}^{\text{loc}}$ and $\widehat{H}_{nk}^{\text{loc}}$ (see Definition 3.1 ahead).

Note that $\widehat{F}_1$ has a jump around $-3$. This jump causes a change of slope in $a_k \widehat{H}_k + a_{K+1} \widehat{H}_+$ for both components $k \in \{1, 2\}$, but only for $k = 1$ is there a touch between $a_k \widehat{H}_k + a_{K+1} \widehat{H}_+$ and $a_k V_k + a_{K+1} V_+$. Similarly, $\widehat{F}_2$ has a jump around $-1$. Again, this causes a change of slope in $a_k \widehat{H}_k + a_{K+1} \widehat{H}_+$ for both components $k \in \{1, 2\}$, but only for $k = 2$ is there a touch between $a_k \widehat{H}_k + a_{K+1} \widehat{H}_+$ and $a_k V_k + a_{K+1} V_+$. The fact that $a_k \widehat{H}_k + a_{K+1} \widehat{H}_+$ has changes of slope without touching $a_k V_k + a_{K+1} V_+$ implies that $a_k \widehat{H}_k + a_{K+1} \widehat{H}_+$ is *not* the convex minorant of $a_k V_k + a_{K+1} V_+$.

It is possible to give convex minorant characterizations of $\widehat{H}$, but again these characterizations are self-induced. Proposition 2.4(a) characterizes $\widehat{H}_k$ in terms of $\sum_{j=1}^{K} \widehat{H}_j$, and Proposition 2.4(b) characterizes $\widehat{H}_k$ in terms of $\sum_{j=1, j \neq k}^{K} \widehat{H}_j$.

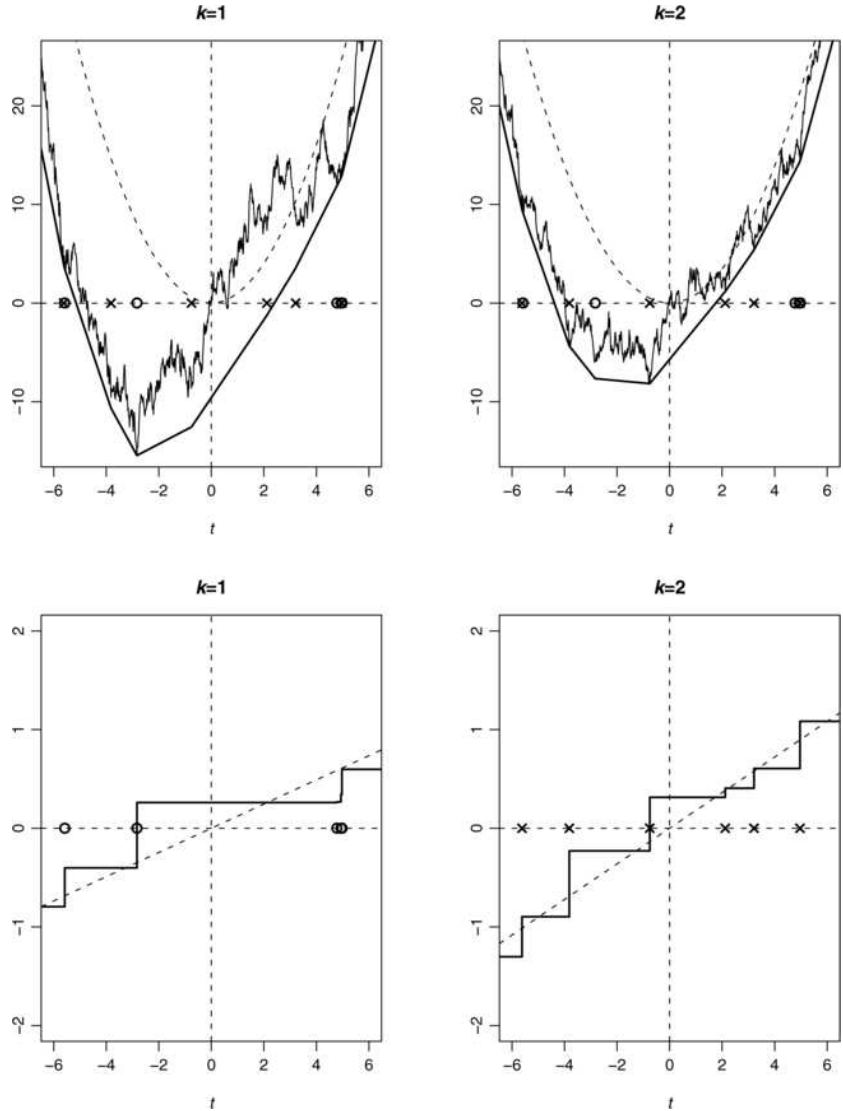

FIG. 1. *Limiting processes for the model given in Example* 2.3 *for* $t_0 = 1$. *The top row shows the processes* $a_k V_k + a_{K+1} V_+$ *and* $a_k \widehat{H}_k + a_{K+1} \widehat{H}_+$, *around the dashed parabolic drifts* $a_k f_{0k}(t_0) t^2/2 + a_{K+1} f_{0+}(t_0) t^2/2$. *The bottom row shows the slope processes* $\widehat{F}_k$, *around dashed lines with slope* $f_{0k}(t_0)$. *The circles and crosses indicate jump points of* $\widehat{F}_1$ *and* $\widehat{F}_2$, *respectively. Note that* $a_k \widehat{H}_k + a_{K+1} \widehat{H}_+$ *touches* $a_k V_k + a_{K+1} V_+$ *whenever* $\widehat{F}_k$ *has a jump, for* $k = 1, 2$.

PROPOSITION 2.4. $\widehat{H}$ satisfies the following convex minorant characterizations:

(a) For each $k = 1, \ldots, K$, $\widehat{H}_k(t)$ is the convex minorant of

$$V_k(t) + \frac{a_{K+1}}{a_k}\{V_+(t) - \widehat{H}_+(t)\}. \tag{8}$$

(b) For each $k = 1, \ldots, K$, $\widehat{H}_k(t)$ is the convex minorant of

$$V_k(t) + \frac{a_{K+1}}{a_k + a_{K+1}}\{V_+^{(-k)}(t) - \widehat{H}_+^{(-k)}(t)\}, \tag{9}$$

where $V_+^{(-k)}(t) = \sum_{j=1, j \neq k}^{K} V_j(t)$ and $\widehat{H}_+^{(-k)}(t) = \sum_{j=1, j \neq k}^{K} \widehat{H}_j(t)$.

PROOF. Conditions (i) and (ii) of Theorem 1.7 are equivalent to

$$\widehat{H}_k(t) \leq V_k(t) + \frac{a_{K+1}}{a_k}\{V_+(t) - \widehat{H}_+(t)\}, \qquad t \in \mathbb{R},$$

$$\int \left\{\widehat{H}_k(t) - V_k(t) - \frac{a_{K+1}}{a_k}\{V_+(t) - \widehat{H}_+(t)\}\right\} d\widehat{F}_k(t) = 0,$$

for $k = 1, \ldots, K$. This gives characterization (a). Similarly, characterization (b) holds since conditions (i) and (ii) of Theorem 1.7 are equivalent to

$$\widehat{H}_k(t) \leq V_k(t) + \frac{a_{K+1}}{a_k + a_{K+1}}\{V_+^{(-k)}(t) - \widehat{H}_+^{(-k)}(t)\}, \qquad t \in \mathbb{R},$$

$$\int \left\{\widehat{H}_k(t) - V_k(t) - \frac{a_{K+1}}{a_k + a_{K+1}}\{V_+^{(-k)}(t) - \widehat{H}_+^{(-k)}(t)\}\right\} d\widehat{F}_k(t) = 0,$$

for $k = 1, \ldots, K$. □

Comparing the MLE and the naive estimator, we see that $\widetilde{H}_k$ is the convex minorant of $V_k$, and $\widehat{H}_k$ is the convex minorant of $V_k + (a_{K+1}/a_k)\{V_+ - \widehat{H}_+\}$. These processes are illustrated in Figure 2. The difference between the two estimators lies in the extra term $(a_{K+1}/a_k)\{V_+ - \widehat{H}_+\}$, which is shown in the bottom row of Figure 2. Apart from the factor $a_{K+1}/a_k$, this term is the same for all $k = 1, \ldots, K$. Furthermore, $a_{K+1}/a_k = F_{0k}(t_0)/F_{0,K+1}(t_0)$ is an increasing function of $t_0$, so that the extra term $(a_{K+1}/a_k)\{V_+ - \widehat{H}_+\}$ is more important for large values of $t_0$. This provides an explanation for the simulation results shown in Figure 3 of Section 4, which indicate that the MLE is superior to the naive estimator in terms of mean squared error, especially for large values of $t$. Finally, note that $(a_{K+1}/a_k)\{V_+ - \widehat{H}_+\}$ appears to be nonnegative in Figure 2. In Proposition 2.5 we prove that this is indeed the case. In turn, this result implies that $\widetilde{H}_k \leq \widehat{H}_k$ (Corollary 2.6), as shown in the top row of Figure 2.

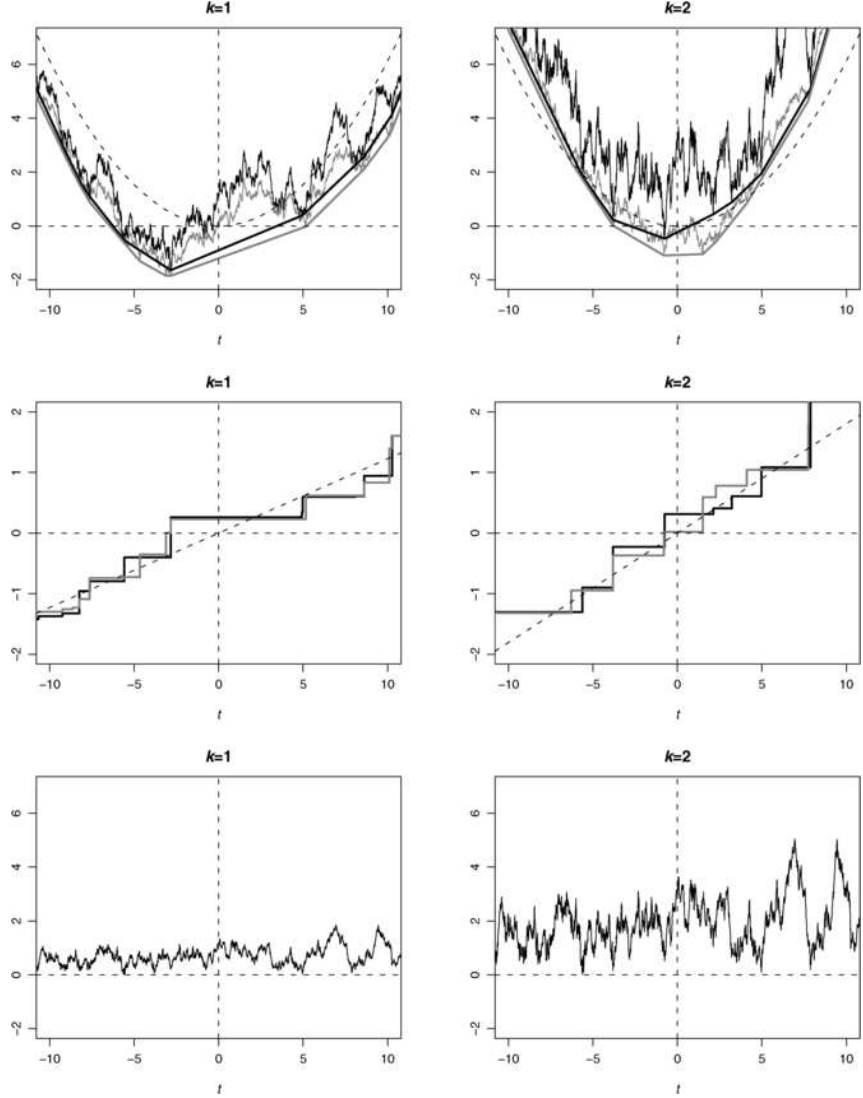

FIG. 2. *Limiting processes for the model given in Example* 2.3 *for* $t_0 = 1$. *The top row shows the processes* $V_k$ *and their convex minorants* $\widetilde{H}_k$ *(gray), together with* $V_k + (a_{K+1}/a_k)(V_+ - \widehat{H}_+)$ *and their convex minorants* $\widehat{H}_k$ *(black). The dashed lines depict the parabolic drift* $f_{0k}(t_0)t^2/2$. *The middle row shows the slope processes* $\widetilde{F}_k$ *(gray) and* $\widehat{F}_k$ *(black), which follow the dashed lines with slope* $f_{0k}(t_0)$. *The bottom row shows the "correction term"* $(a_{K+1}/a_k)(V_+ - \widehat{H}_+)$ *for the MLE.*

PROPOSITION 2.5. $\widehat{H}_+(t) \leq V_+(t)$ for all $t \in \mathbb{R}$.

PROOF. Theorem 1.7(i) can be written as

$$\widehat{H}_k(t) + \frac{F_{0k}(t_0)}{1 - F_{0+}(t_0)}\widehat{H}_+(t) \leq V_k(t) + \frac{F_{0k}(t_0)}{1 - F_{0+}(t_0)}V_+(t)$$

$$\text{for } k = 1, \ldots, K, t \in \mathbb{R}.$$

The statement then follows by summing over $k = 1, \ldots, K$. □

COROLLARY 2.6. $\widetilde{H}_k(t) \leq \widehat{H}_k(t)$ for all $k = 1, \ldots, K$ and $t \in \mathbb{R}$.

PROOF. Let $k \in \{1, \ldots, K\}$ and recall that $\widetilde{H}_k$ is the convex minorant of $V_k$. Since $V_+ - \widehat{H}_+ \geq 0$ by Proposition 2.5, it follows that $\widetilde{H}_k$ is a convex function below $V_k + (a_{K+1}/a_k)\{V_+ - \widehat{H}_+\}$. Hence, it is bounded above by the convex minorant $\widehat{H}_k$ of $V_k + (a_{K+1}/a_k)\{V_+ - \widehat{H}_+\}$. □

Finally, we write the characterization of Theorem 1.7 in a way that is analogous to the characterization of the MLE in Proposition 4.8 of [8]. We do this to make a connection between the finite sample situation and the limiting situation. Using this connection, the proofs for the tightness results in Section 2.2 are similar to the proofs for the local rate of convergence in [8], Section 4.3. We need the following definition:

DEFINITION 2.7. For $k = 1, \ldots, K$ and $t \in \mathbb{R}$, we define

(10) $\quad \bar{F}_{0k}(t) = f_{0k}(t_0)t \quad \text{and} \quad S_k(t) = a_k W_k(t) + a_{K+1} W_+(t).$

Note that $S_k$ is the limit of a rescaled version of the process $S_{nk} = a_k W_{nk} + a_{K+1} W_{n+}$, defined in (18) of [8].

PROPOSITION 2.8. For all $k = 1, \ldots, K$, for each point $\tau_k \in \widehat{\mathcal{N}}_k$ [defined in (7)] and for all $s \in \mathbb{R}$, we have

(11) $\quad \int_{\tau_k}^{s} \{a_k\{\widehat{F}_k(u) - \bar{F}_{0k}(u)\} + a_{K+1}\{\widehat{F}_+(u) - \bar{F}_{0+}(u)\}\} \, du \leq \int_{\tau_k}^{s} dS_k(u),$

and equality must hold if $s \in \widehat{\mathcal{N}}_k$.

PROOF. Let $k \in \{1, \ldots, K\}$. By Theorem 1.7(i), we have

$$a_k \widehat{H}_k(t) + a_{K+1}\widehat{H}_+(t) \leq a_k V_k(t) + a_{K+1}V_+(t), \qquad t \in \mathbb{R},$$

where equality holds at $t = \tau_k \in \widehat{\mathcal{N}}_k$. Subtracting this expression for $t = \tau_k$ from the expression for $t = s$, we get

$$\int_{\tau_k}^{s} \{a_k \widehat{F}_k(u) + a_{K+1} \widehat{F}_+(u)\} \, du \leq \int_{\tau_k}^{s} \{a_k \, dV_k(u) + a_{K+1} \, dV_+(u)\}.$$

The result then follows by subtracting $\int_{\tau_k}^{s} \{a_k \bar{F}_{0k}(u) + a_{K+1} \bar{F}_{0+}(u)\} \, du$ from both sides, and using that $dV_k(u) = \bar{F}_{0k}(u) \, du + dW_k(u)$ [see (4)]. □

2.2. *Tightness of $\widehat{H}$ and $\widehat{F}$.* The main results of this section are tightness of $\{\widehat{F}_k(t) - \bar{F}_{0k}(t)\}$ (Proposition 2.9) and $\{\widehat{H}_k(t) - V_k(t)\}$ (Corollary 2.15), for $t \in \mathbb{R}$. These results are used in Section 2.3 to prove that $\widehat{H}$ and $\widehat{F}$ are almost surely unique.

PROPOSITION 2.9. *For every $\varepsilon > 0$ there is an $M > 0$ such that*

$$P(|\widehat{F}_k(t) - \bar{F}_{0k}(t)| \geq M) < \varepsilon \qquad \textit{for } k = 1, \ldots, K, t \in \mathbb{R}.$$

Proposition 2.9 is the limit version of Theorem 4.17 of [8], which gave the $n^{1/3}$ local rate of convergence of $\widehat{F}_{nk}$. Hence, analogously to [8], proof of Theorem 4.17, we first prove a stronger tightness result for the sum process $\{\widehat{F}_+(t) - \bar{F}_{0+}(t)\}$, $t \in \mathbb{R}$.

PROPOSITION 2.10. *Let $\beta \in (0, 1)$ and define*

$$(12) \qquad v(t) = \begin{cases} 1, & \textit{if } |t| \leq 1, \\ |t|^\beta, & \textit{if } |t| > 1. \end{cases}$$

*Then for every $\varepsilon > 0$ there is an $M > 0$ such that*

$$P\left(\sup_{t \in \mathbb{R}} \frac{|\widehat{F}_+(t) - \bar{F}_{0+}(t)|}{v(t - s)} \geq M\right) < \varepsilon \qquad \textit{for } s \in \mathbb{R}.$$

PROOF. The organization of this proof is similar to the proof of Theorem 4.10 of [8]. Let $\varepsilon > 0$. We only prove the result for $s = 0$, since the proof for $s \neq 0$ is equivalent, due to stationarity of the increments of Brownian motion.

It is sufficient to show that we can choose $M > 0$ such that

$$P(\exists t \in \mathbb{R} : \widehat{F}_+(t) \notin (\bar{F}_{0+}(t - Mv(t)), \bar{F}_{0+}(t + Mv(t))))$$
$$= P(\exists t \in \mathbb{R} : |\widehat{F}_+(t) - \bar{F}_{0+}(t)| \geq f_{0+}(t_0) Mv(t)) < \varepsilon.$$

In fact, we only prove that there is an $M$ such that

$$P(\exists t \in [0, \infty) : \widehat{F}_+(t) \geq \bar{F}_{0+}(t + Mv(t))) < \frac{\varepsilon}{4},$$

since the proofs for the inequality $\widehat{F}_+(t) \leq \bar{F}_{0+}(t - Mv(t))$ and the interval $(-\infty, 0]$ are analogous. In turn, it is sufficient to show that there is an $m_1 > 0$ such that

(13) $P(\exists t \in [j, j+1) : \widehat{F}_+(t) \geq \bar{F}_{0+}(t + Mv(t))) \leq p_{jM}, \qquad j \in \mathbb{N}, M > m_1,$

where $p_{jM}$ satisfies $\sum_{j=0}^{\infty} p_{jM} \to 0$ as $M \to \infty$. We prove (13) for

(14) $$p_{jM} = d_1 \exp\{-d_2(Mv(j))^3\},$$

where $d_1$ and $d_2$ are positive constants. Using the monotonicity of $\widehat{F}_+$, we only need to show that $P(A_{jM}) \leq p_{jM}$ for all $j \in \mathbb{N}$ and $M > m_1$, where

(15) $\quad A_{jM} = \{\widehat{F}_+(j+1) \geq \bar{F}_{0+}(s_{jM})\} \quad \text{and} \quad s_{jM} = j + Mv(j).$

We now fix $M > 0$ and $j \in \mathbb{N}$, and define $\tau_{kj} = \max\{\widehat{\mathcal{N}}_k \cap (-\infty, j+1]\}$, for $k = 1, \ldots, K$. These points are well defined by Theorem 1.7(iii) and Corollary 2.2(i). Without loss of generality, we assume that the sub-distribution functions are labeled so that $\tau_{1j} \leq \cdots \leq \tau_{Kj}$. On the event $A_{jM}$, there is a $k \in \{1, \ldots, K\}$ such that $\widehat{F}_k(j+1) \geq \bar{F}_{0k}(s_{jM})$. Hence, we can define $\ell \in \{1, \ldots, K\}$ such that

(16) $\qquad \widehat{F}_k(j+1) < \bar{F}_{0k}(s_{jM}), \qquad k = \ell+1, \ldots, K,$

(17) $\qquad \widehat{F}_\ell(j+1) \geq \bar{F}_{0\ell}(s_{jM}).$

Recall that $\widehat{F}$ must satisfy (11). Hence, $P(A_{jM})$ equals

$$P\biggl(\int_{\tau_{\ell j}}^{s_{jM}} \{a_\ell\{\widehat{F}_\ell(u) - \bar{F}_{0\ell}(u)\} + a_{K+1}\{\widehat{F}_+(u) - \bar{F}_{0+}(u)\}\} \, du$$
$$\leq \int_{\tau_{\ell j}}^{s_{jM}} dS_\ell(u), A_{jM}\biggr)$$

(18) $\qquad \leq P\biggl(\int_{\tau_{\ell j}}^{s_{jM}} a_\ell\{\widehat{F}_\ell(u) - \bar{F}_{0\ell}(u)\} \, du \leq \int_{\tau_{\ell j}}^{s_{jM}} dS_\ell(u), A_{jM}\biggr)$

(19) $\qquad + P\biggl(\int_{\tau_{\ell j}}^{s_{jM}} \{\widehat{F}_+(u) - \bar{F}_{0+}(u)\} \, du \leq 0, A_{jM}\biggr).$

Using the definition of $\tau_{\ell j}$ and the fact that $\widehat{F}_\ell$ is monotone nondecreasing and piecewise constant (Corollary 2.2), it follows that on the event $A_{jM}$ we have $\widehat{F}_\ell(u) \geq \widehat{F}_\ell(\tau_{\ell j}) = \widehat{F}_\ell(j+1) \geq \bar{F}_{0\ell}(s_{jM})$, for $u \geq \tau_{\ell j}$. Hence, we can bound (18) above by

$$P\biggl(\int_{\tau_{\ell j}}^{s_{jM}} a_\ell\{\bar{F}_{0\ell}(s_{jM}) - \bar{F}_{0\ell}(u)\} \, du \leq \int_{\tau_{\ell j}}^{s_{jM}} dS_\ell(u)\biggr)$$

$$= P\bigg(\tfrac{1}{2}f_{0\ell}(t_0)(s_{jM} - \tau_{\ell j})^2 \leq \int_{\tau_{\ell j}}^{s_{jM}} dS_\ell(u)\bigg)$$

$$\leq P\bigg(\inf_{w \leq j+1}\bigg\{\tfrac{1}{2}f_{0\ell}(t_0)(s_{jM} - w)^2 - \int_w^{s_{jM}} dS_\ell(u)\bigg\} \leq 0\bigg).$$

For $m_1$ sufficiently large, this probability is bounded above by $p_{jM}/2$ for all $M > m_1$ and $j \in \mathbb{N}$, by Lemma 2.11 below. Similarly, (19) is bounded by $p_{jM}/2$, using Lemma 2.12 below. $\square$

Lemmas 2.11 and 2.12 are the key lemmas in the proof of Proposition 2.10. They are the limit versions of Lemmas 4.13 and 4.14 of [8], and their proofs are given in Section 5. The basic idea of Lemma 2.11 is that the positive quadratic drift $b(s_{jM} - w)^2$ dominates the Brownian motion process $S_k$ and the term $C(s_{jM} - w)^{3/2}$. Note that the lemma also holds when $C(s_{jM} - w)^{3/2}$ is omitted, since this term is positive for $M > 1$. In fact, in the proof of Proposition 2.10 we only use the lemma without this term, but we need the term $C(s_{jM} - w)^{3/2}$ in the proof of Proposition 2.9 ahead. The proof of Lemma 2.12 relies on the system of component processes. Since it is very similar to the proof of Lemma 4.14, we only point out the differences in Section 5.

LEMMA 2.11. *Let $C > 0$ and $b > 0$. Then there exists an $m_1 > 0$ such that for all $k = 1, \ldots, K$, $M > m_1$ and $j \in \mathbb{N}$,*

$$P\bigg(\inf_{w \leq j+1}\bigg\{b(s_{jM} - w)^2 - \int_w^{s_{jM}} dS_k(u) - C(s_{jM} - w)^{3/2}\bigg\} \leq 0\bigg) \leq p_{jM},$$

*where $s_{jM} = j + Mv(j)$, and $S_k(\cdot)$, $v(\cdot)$ and $p_{jM}$ are defined by (10), (12) and (14), respectively.*

LEMMA 2.12. *Let $\ell$ be defined by (16) and (17). There is an $m_1 > 0$ such that*

$$P\bigg(\int_{\tau_{\ell j}}^{s_{jM}}\{\widehat{F}_+(u) - \bar{F}_{0+}(u)\}\,du \leq 0, A_{jM}\bigg) \leq p_{jM} \qquad \text{for } M > m_1, j \in \mathbb{N},$$

*where $s_{jM} = j + Mv(j)$, $\tau_{\ell j} = \max\{\widehat{\mathcal{N}}_\ell \cap (-\infty, j+1]\}$, and $v(\cdot)$, $p_{jM}$ and $A_{jM}$ are defined by (12), (14) and (15), respectively.*

In order to prove tightness of $\{\widehat{F}_k(t) - \bar{F}_{0k}(t)\}$, $t \in \mathbb{R}$, we only need Proposition 2.10 to hold for one value of $\beta \in (0, 1)$, analogously to [8], Remark 4.12. We therefore fix $\beta = 1/2$, so that $v(t) = 1 \vee \sqrt{|t|}$. Then Proposition 2.10 leads to the following corollary, which is a limit version of Corollary 4.16 of [8]:

COROLLARY 2.13. *For every $\varepsilon > 0$ there is a $C > 0$ such that*

$$P\left\{\sup_{u \in \mathbb{R}_+} \frac{\int_{s-u}^{s} |\widehat{F}_+(t) - \bar{F}_{0+}(t)|\, dt}{u \vee u^{3/2}} \geq C\right\} < \varepsilon \qquad \text{for } s \in \mathbb{R}.$$

This corollary allows us to complete the proof of Proposition 2.9.

PROOF OF PROPOSITION 2.9. Let $\varepsilon > 0$ and let $k \in \{1, \ldots, K\}$. It is sufficient to show that there is an $M > 0$ such that $P(\widehat{F}_k(t) \geq \bar{F}_{0k}(t+M)) < \varepsilon$ and $P(\widehat{F}_k(t) \leq \bar{F}_{0k}(t-M)) < \varepsilon$ for all $t \in \mathbb{R}$. We only prove the first inequality, since the proof of the second one is analogous. Thus, let $t \in \mathbb{R}$ and $M > 1$, and define

$$B_{kM} = \{\widehat{F}_k(t) \geq \bar{F}_{0k}(t+M)\} \quad \text{and} \quad \tau_k = \max\{\widehat{\mathcal{N}}_k \cap (-\infty, t]\}.$$

Note that $\tau_k$ is well defined because of Theorem 1.7(iii) and Corollary 2.2(i). We want to prove that $P(B_{kM}) < \varepsilon$. Recall that $\widehat{F}$ must satisfy (11). Hence,

$$(20) \quad P(B_{kM}) = P\left(\int_{\tau_k}^{t+M} \{a_k\{\widehat{F}_k(u) - \bar{F}_{0k}(u)\} + a_{K+1}\{\widehat{F}_+(u) - \bar{F}_{0+}(u)\}\}\, du\right.$$
$$\left. \leq \int_{\tau_k}^{t+M} dS_k(u), B_{kM}\right).$$

By Corollary 2.13, we can choose $C > 0$ such that, with high probability,

$$(21) \quad \int_{\tau_k}^{t+M} |\widehat{F}_+(u) - \bar{F}_{0+}(u)|\, du \leq C(t + M - \tau_k)^{3/2},$$

uniformly in $\tau_k \leq t$, using that $u^{3/2} > u$ for $u > 1$. Moreover, on the event $B_{kM}$, we have $\int_{\tau_k}^{t+M}\{\widehat{F}_k(u) - \bar{F}_{0k}(u)\}\, du \geq \int_{\tau_k}^{t+M}\{\bar{F}_{0k}(t+M) - \bar{F}_{0k}(u)\}\, du = f_{0k}(t_0)(t+M-\tau_k)^2/2$, yielding a positive quadratic drift. The statement now follows by combining these facts with (20), and applying Lemma 2.11. □

Proposition 2.9 leads to the following corollary about the distance between the jump points of $\widehat{F}_k$. The proof is analogous to the proof of Corollary 4.19 of [8], and is therefore omitted.

COROLLARY 2.14. *For all $k = 1, \ldots, K$, let $\tau_k^-(s)$ and $\tau_k^+(s)$ be, respectively, the largest jump point $\leq s$ and the smallest jump point $> s$ of $\widehat{F}_k$. Then for every $\varepsilon > 0$ there is a $C > 0$ such that $P(\tau_k^+(s) - \tau_k^-(s) > C) < \varepsilon$, for $k = 1, \ldots, K$, $s \in \mathbb{R}$.*

Combining Theorem 2.9 and Corollary 2.14 yields tightness of $\{\widehat{H}_k(t) - V_k(t)\}$:

COROLLARY 2.15. *For every $\varepsilon > 0$ there is an $M > 0$ such that*

$$P(|\widehat{H}_k(t) - V_k(t)| > M) < \varepsilon \qquad \text{for } t \in \mathbb{R}.$$

2.3. *Uniqueness of $\widehat{H}$ and $\widehat{F}$.* We now use the tightness results of Section 2.2 to prove the uniqueness part of Theorem 1.7, as given in Proposition 2.16. The existence part of Theorem 1.7 will follow in Section 3.

PROPOSITION 2.16. *Let $\widehat{H}$ and $H$ satisfy the conditions of Theorem 1.7. Then $\widehat{H} \equiv H$ almost surely.*

The proof of Proposition 2.16 relies on the following lemma:

LEMMA 2.17. *Let $\widehat{H} = (\widehat{H}_1, \ldots, \widehat{H}_K)$ and $H = (H_1, \ldots, H_K)$ satisfy the conditions of Theorem 1.7, and let $\widehat{F} = (\widehat{F}_1, \ldots, \widehat{F}_K)$ and $F = (F_1, \ldots, F_K)$ be the corresponding derivatives. Then*

$$\begin{aligned}
(22) \quad & \sum_{k=1}^{K} a_k \int \{F_k(t) - \widehat{F}_k(t)\}^2 \, dt + a_{K+1} \int \{F_+(t) - \widehat{F}_+(t)\}^2 \, dt \\
& \leq \liminf_{m \to \infty} \sum_{k=1}^{K} \{\psi_k(m) - \psi_k(-m)\},
\end{aligned}$$

*where $\psi_k : \mathbb{R} \to \mathbb{R}$ is defined by*

$$(23) \quad \psi_k(t) = \{F_k(t) - \widehat{F}_k(t)\}[a_k\{H_k(t) - \widehat{H}_k(t)\} + a_{K+1}\{H_+(t) - \widehat{H}_+(t)\}].$$

PROOF. We define the following functional:

$$\begin{aligned}
\phi_m(F) = & \sum_{k=1}^{K} a_k \left\{ \tfrac{1}{2} \int_{-m}^{m} F_k^2(t) \, dt - \int_{-m}^{m} F_k(t) \, dV_k(t) \right\} \\
& + a_{K+1} \left\{ \tfrac{1}{2} \int_{-m}^{m} F_+^2(t) \, dt - \int_{-m}^{m} F_+(t) \, dV_+(t) \right\}, \qquad m \in \mathbb{N}.
\end{aligned}$$

Then, letting

$$(24) \qquad D_k(t) = a_k\{H_k(t) - V_k(t)\} + a_{K+1}\{H_+(t) - V_+(t)\},$$

$$(25) \qquad \widehat{D}_k(t) = a_k\{\widehat{H}_k(t) - V_k(t)\} + a_{K+1}\{\widehat{H}_+(t) - V_+(t)\},$$

and using $F_k^2 - \widehat{F}_k^2 = (F_k - \widehat{F}_k)^2 + 2\widehat{F}_k(F_k - \widehat{F}_k)$, we have

$$\phi_m(F) - \phi_m(\widehat{F}) = \sum_{k=1}^{K} \frac{a_k}{2} \int_{-m}^{m} \{F_k(t) - \widehat{F}_k(t)\}^2 \, dt$$

(26)
$$+ \frac{a_{K+1}}{2} \int_{-m}^{m} \{F_+(t) - \widehat{F}_+(t)\}^2 \, dt$$

$$+ \sum_{k=1}^{K} \int_{-m}^{m} \{F_k(t) - \widehat{F}_k(t)\} \, d\widehat{D}_k(t).$$

Using integration by parts, we rewrite the last term of the right-hand side of (26) as:

(27)
$$\sum_{k=1}^{K} \{F_k(t) - \widehat{F}_k(t)\}\widehat{D}_k(t)\Big|_{-m}^{m} - \sum_{k=1}^{K} \int_{-m}^{m} \widehat{D}_k(t) \, d\{F_k(t) - \widehat{F}_k(t)\}$$

$$\geq \sum_{k=1}^{K} \{F_k(t) - \widehat{F}_k(t)\}\widehat{D}_k(t)\Big|_{-m}^{m}.$$

The inequality on the last line Follows from: (a) $\int_{-m}^{m} \widehat{D}_k(t) \, d\widehat{F}_k(t) = 0$ by Theorem 1.7(ii), and (b) $\int_{-m}^{m} \widehat{D}_k(t) \, dF_k(t) \leq 0$, since $\widehat{D}_k(t) \leq 0$ by Theorem 1.7(i) and $F_k$ is monotone nondecreasing. Combining (26) and (27), and using the same expressions with $F$ and $\widehat{F}$ interchanged, yields

$$0 = \phi_m(\widehat{F}) - \phi_m(F) + \phi_m(F) - \phi_m(\widehat{F})$$

$$\geq \sum_{k=1}^{K} a_k \int_{-m}^{m} \{F_k(t) - \widehat{F}_k(t)\}^2 \, dt + a_{K+1} \int_{-m}^{m} \{F_+(t) - \widehat{F}_+(t)\}^2 \, dt$$

$$+ \sum_{k=1}^{K} \{\widehat{F}_k(t) - F_k(t)\}D_k(t)\Big|_{-m}^{m} + \sum_{k=1}^{K} \{F_k(t) - \widehat{F}_k(t)\}\widehat{D}_k(t)\Big|_{-m}^{m}.$$

By writing out the right-hand side of this expression, we find that it is equivalent to

$$\sum_{k=1}^{K} a_k \int_{-m}^{m} \{F_k(t) - \widehat{F}_k(t)\}^2 \, dt + a_{K+1} \int_{-m}^{m} \{F_+(t) - \widehat{F}_+(t)\}^2 \, dt$$

(28)
$$\leq \sum_{k=1}^{K} [\{F_k(m) - \widehat{F}_k(m)\}\{D_k(m) - \widehat{D}_k(m)\}$$

$$- \{F_k(-m) - \widehat{F}_k(-m)\}\{D_k(-m) - \widehat{D}_k(-m)\}].$$

This inequality holds for all $m \in \mathbb{N}$, and hence we can take $\liminf_{m \to \infty}$. The left-hand side of (28) is a monotone sequence in $m$, so that we can replace $\liminf_{m\to\infty}$ by $\lim_{m\to\infty}$. The result then follows from the definitions of $\psi_k$, $D_k$ and $\widehat{D}_k$ in (23)–(25). □

We are now ready to prove Proposition 2.16. The idea of the proof is to show that the right-hand side of (22) is almost surely equal to zero. We prove this in two steps. First, we show that it is of order $O_p(1)$, using the tightness results of Proposition 2.9 and Corollary 2.15. Next, we show that the right-hand side is almost surely equal to zero.

PROOF OF PROPOSITION 2.16. We first show that the right-hand side of (22) is of order $O_p(1)$. Let $k \in \{1,\ldots,K\}$, and note that Proposition 2.9 yields that $\{F_k(m) - \bar{F}_{0k}(m)\}$ and $\{\widehat{F}_k(m) - \bar{F}_{0k}(m)\}$ are of order $O_p(1)$, so that also $\{F_k(m) - \widehat{F}_k(m)\} = O_p(1)$. Similarly, Corollary 2.15 implies that $\{H_k(m) - \widehat{H}_k(m)\} = O_p(1)$. Using the same argument for $-m$, this proves that the right-hand side of (22) is of order $O_p(1)$.

We now show that the right-hand side of (22) is almost surely equal to zero. Let $k \in \{1,\ldots,K\}$. We only consider $|F_k(m) - \widehat{F}_k(m)||H_k(m) - \widehat{H}_k(m)|$, since the term $|F_k(m) - \widehat{F}_k(m)||H_+(m) - \widehat{H}_+(m)|$ and the point $-m$ can be treated analogously. It is sufficient to show that

$$(29) \quad \liminf_{m\to\infty} P(|F_k(m) - \widehat{F}_k(m)||H_k(m) - \widehat{H}_k(m)| > \eta) = 0 \quad \text{for all } \eta > 0.$$

Let $\tau_{mk}$ be the last jump point of $F_k$ before $m$, and let $\widehat{\tau}_{mk}$ be the last jump point of $\widehat{F}_k$ before $m$. We define the following events:

$$E_m = E_m(\varepsilon, \delta, C) = E_{1m}(\varepsilon) \cap E_{2m}(\delta) \cap E_{3m}(C) \quad \text{where}$$

$$E_{1m} = E_{1m}(\varepsilon) = \left\{ \int_{\tau_{mk}\vee\widehat{\tau}_{mk}}^{\infty} \{F_k(t) - \widehat{F}_k(t)\}^2 \, dt < \varepsilon \right\},$$

$$E_{2m} = E_{2m}(\delta) = \{m - (\tau_{mk} \vee \widehat{\tau}_{mk}) > \delta\},$$

$$E_{3m} = E_{3m}(C) = \{|H_k(m) - \widehat{H}_k(m)| < C\}.$$

Let $\varepsilon_1 > 0$ and $\varepsilon_2 > 0$. Since the right-hand side of (22) is of order $O_p(1)$, it follows that $\int \{F_k(t) - \widehat{F}_k(t)\}^2 \, dt = O_p(1)$ for every $k \in \{1,\ldots,K\}$. This implies that $\int_m^{\infty} \{F_k(t) - \widehat{F}_k(t)\}^2 \, dt \to_p 0$ as $m \to \infty$. Together with the fact that $m - \{\tau_{mk} \vee \widehat{\tau}_{mk}\} = O_p(1)$ (Corollary 2.14), this implies that there is an $m_1 > 0$ such that $P(E_{1m}(\varepsilon_1)^c) < \varepsilon_1$ for all $m > m_1$. Next, recall that the points of jump of $F_k$ and $\widehat{F}_k$ are contained in the set $\mathcal{N}_k$, defined in Proposition 2.1. Letting $\tau'_{mk} = \max\{\mathcal{N}_k \cap (-\infty, m]\}$, we have

$$(30) \quad P(E_{2m}^c(\delta)) \leq P(m - \tau'_{mk} < \delta).$$

The distribution of $m - \tau'_{mk}$ is independent of $m$, nondegenerate and continuous (see [4]). Hence, we can choose $\delta > 0$ such that the probabilities in

(30) are bounded by $\varepsilon_2/2$ for all $m$. Furthermore, by tightness of $\{H_k(m) - \widehat{H}_k(m)\}$, there is a $C > 0$ such that $P(E_{3m}(C)^c) < \varepsilon_2/2$ for all $m$. This implies that $P(E_m(\varepsilon_1, \delta, C)^c) < \varepsilon_1 + \varepsilon_2$ for $m > m_1$.

Returning to (29), we now have for $\eta > 0$:

$$\liminf_{m \to \infty} P(|F_k(m) - \widehat{F}_k(m)||H_k(m) - \widehat{H}_k(m)| > \eta)$$
$$\leq \varepsilon_1 + \varepsilon_2$$
$$\quad + \liminf_{m \to \infty} P(|F_k(m) - \widehat{F}_k(m)||H_k(m) - \widehat{H}_k(m)| > \eta, E_m(\varepsilon_1, \delta, C))$$
$$\leq \varepsilon_1 + \varepsilon_2 + \liminf_{m \to \infty} P\left(|F_k(m) - \widehat{F}_k(m)| > \frac{\eta}{C}, E_m(\varepsilon_1, \delta, C)\right),$$

using the definition of $E_{3m}(C)$ in the last line. The probability in the last line equals zero for $\varepsilon_1$ small. To see this, note that $F_k(m) - \widehat{F}_k(m) > \eta/C$, $m - \{\tau_{mk} \vee \widehat{\tau}_{mk}\} > \delta$, and the fact that $F_k$ and $\widehat{F}_k$ are piecewise constant on $m - \{\tau_{km} \vee \widehat{\tau}_{km}\}$ imply that

$$\int_{\tau_{mk} \vee \widehat{\tau}_{mk}}^{\infty} \{F_k(u) - \widehat{F}_k(u)\}^2 \, du \geq \int_{\tau_{mk} \vee \widehat{\tau}_{mk}}^{m} \{F_k(u) - \widehat{F}_k(u)\}^2 \, du > \frac{\eta^2 \delta}{C^2},$$

so that $E_{1m}(\varepsilon_1)$ cannot hold for $\varepsilon_1 < \eta^2 \delta/C^2$.

This proves that the right-hand side of (22) equals zero, almost surely. Together with the right-continuity of $F_k$ and $\widehat{F}_k$, this implies that $F_k \equiv \widehat{F}_k$ almost surely, for $k = 1, \ldots, K$. Since $F_k$ and $\widehat{F}_k$ are the right derivatives of $H_k$ and $\widehat{H}_k$, this yields that $H_k = \widehat{H}_k + c_k$ almost surely. Finally, both $H_k$ and $\widehat{H}_k$ satisfy conditions (i) and (ii) of Theorem 1.7 for $k = 1, \ldots, K$, so that $c_1 = \cdots = c_K = 0$ and $H \equiv \widehat{H}$ almost surely. □

**3. Proof of the limiting distribution of the MLE.** In this section we prove that the MLE converges to the limiting distribution given in Theorem 1.8. In the process, we also prove the existence part of Theorem 1.7.

First, we recall from [8], Section 2.2, that the naive estimators $\widetilde{F}_{nk}$, $k = 1, \ldots, K$, are unique at $t \in \{T_1, \ldots, T_n\}$, and that the MLEs $\widehat{F}_{nk}$, $k = 1, \ldots, K$, are unique at $t \in \mathcal{T}_K$, where $\mathcal{T}_k = \{T_i, i = 1, \ldots, n : \Delta_k^i + \Delta_{K+1}^i > 0\} \cup \{T_{(n)}\}$ for $k = 1, \ldots, K$ (see [8], Proposition 2.3). To avoid issues with non-uniqueness, we adopt the convention that $\widetilde{F}_{nk}$ and $\widehat{F}_{nk}$, $k = 1, \ldots, K$, are piecewise constant and right-continuous, with jumps only at the points at which they are uniquely defined. This convention does not affect the asymptotic properties of the estimators under the assumptions of Section 1.2. Recalling the definitions of $G$ and $G_n$ given in Section 1.1, we now define the following localized processes:

DEFINITION 3.1. For each $k = 1, \ldots, K$, we define

$$\widehat{F}_{nk}^{\mathrm{loc}}(t) = n^{1/3}\{\widehat{F}_{nk}(t_0 + n^{-1/3}t) - F_{0k}(t_0)\}, \tag{31}$$

$$V_{nk}^{\mathrm{loc}}(t) = \frac{n^{2/3}}{g(t_0)} \int_{u \in (t_0, t_0 + n^{-1/3}t]} \{\delta_k - F_{0k}(t_0)\} \, d\mathbb{P}_n(u, \delta), \tag{32}$$

$$\bar{H}_{nk}^{\mathrm{loc}}(t) = \frac{n^{2/3}}{g(t_0)} \int_{t_0}^{t_0 + n^{-1/3}t} \{\widehat{F}_{nk}(u) - F_{0k}(t_0)\} \, dG(u), \tag{33}$$

$$\widehat{H}_{nk}^{\mathrm{loc}}(t) = \bar{H}_{nk}^{\mathrm{loc}}(t) + \frac{c_{nk}}{a_k} - F_{0k}(t_0) \sum_{k=1}^{K} \frac{c_{nk}}{a_k}, \tag{34}$$

where $c_{nk}$ is the difference between $a_k V_{nk}^{\mathrm{loc}} + a_{K+1} V_{n+}^{\mathrm{loc}}$ and $a_k H_{nk}^{\mathrm{loc}} + a_{K+1} H_{n+}^{\mathrm{loc}}$ at the last jump point $\tau_{nk}$ of $\widehat{F}_{nk}^{\mathrm{loc}}$ before zero, that is,

$$c_{nk} = a_k V_{nk}^{\mathrm{loc}}(\tau_{nk}-) + a_{K+1} V_{n+}^{\mathrm{loc}}(\tau_{nk}-) - a_k \bar{H}_{nk}^{\mathrm{loc}}(\tau_{nk}) - a_{K+1} \bar{H}_{n+}^{\mathrm{loc}}(\tau_{nk}). \tag{35}$$

Moreover, we define the vectors $\widehat{F}_n^{\mathrm{loc}} = (\widehat{F}_{n1}^{\mathrm{loc}}, \ldots, \widehat{F}_{nK}^{\mathrm{loc}})$, $V_n^{\mathrm{loc}} = (V_{n1}^{\mathrm{loc}}, \ldots, V_{nK}^{\mathrm{loc}})$ and $\widehat{H}_n^{\mathrm{loc}} = (\widehat{H}_{n1}^{\mathrm{loc}}, \ldots, \widehat{H}_{nK}^{\mathrm{loc}})$.

Note that $\widehat{H}_{nk}^{\mathrm{loc}}$ differs from $\bar{H}_{nk}^{\mathrm{loc}}$ by a vertical shift, and that $(\widehat{H}_{nk}^{\mathrm{loc}})'(t) = (\bar{H}_{nk}^{\mathrm{loc}})'(t) = \widehat{F}_{nk}^{\mathrm{loc}}(t) + o(1)$. We now show that the MLE satisfies the characterization given in Proposition 3.2, which can be viewed as a recentered and rescaled version of the characterization in Proposition 4.8 of [8]. In the proof of Theorem 1.8 we will see that, as $n \to \infty$, this characterization converges to the characterization of the limiting process given in Theorem 1.7.

PROPOSITION 3.2. *Let the assumptions of Section* 1.2 *hold, and let* $m > 0$. *Then*

$$a_k \widehat{H}_{nk}^{\mathrm{loc}}(t) + a_{K+1} \widehat{H}_{n+}^{\mathrm{loc}}(t)$$
$$\leq a_k V_{nk}^{\mathrm{loc}}(t-) + a_{K+1} V_{n+}^{\mathrm{loc}}(t-) + R_{nk}^{\mathrm{loc}}(t) \quad \text{for } t \in [-m, m],$$

$$\int_{-m}^{m} \{a_k V_{nk}^{\mathrm{loc}}(t-) + a_{K+1} \, dV_{n+}^{\mathrm{loc}}(t-)$$
$$+ R_{nk}^{\mathrm{loc}}(t) - a_k \widehat{H}_{nk}^{\mathrm{loc}}(t) - a_{K+1} \widehat{H}_{n+}^{\mathrm{loc}}(t)\} \, d\widehat{F}_{nk}^{\mathrm{loc}}(t) = 0,$$

*where* $R_{nk}^{\mathrm{loc}}(t) = o_p(1)$, *uniformly in* $t \in [-m, m]$.

PROOF. Let $m > 0$ and let $\tau_{nk}$ be the last jump point of $\widehat{F}_{nk}$ before $t_0$. It follows from the characterization of the MLE in Proposition 4.8 of [8] that

$$\int_{\tau_{nk}}^{s} \{a_k\{\widehat{F}_{nk}(u) - F_{0k}(u)\} + a_{K+1}\{\widehat{F}_{n+}(u) - F_{0+}(u)\}\} \, dG(u)$$

$$(36) \qquad \leq \int_{[\tau_{nk},s)} \{a_k\{\delta_k - F_{0k}(u)\} + a_{K+1}\{\delta_+ - F_{0+}(u)\}\} d\mathbb{P}_n(u,\delta)$$
$$+ R_{nk}(\tau_{nk}, s),$$

where equality holds if $s$ is a jump point of $\widehat{F}_{nk}$. Using that $t_0 - \tau_{nk} = O_p(n^{-1/3})$ by [8], Corollary 4.19, it follows from [8], Corollary 4.20 that $R_{nk}(\tau_{nk}, s) = o_p(n^{-2/3})$, uniformly in $s \in [t_0 - m_1 n^{1/3}, t_0 + m_1 n^{-1/3}]$. We now add

$$\int_{\tau_{nk}}^{s} \{a_k\{F_{0k}(u) - F_{0k}(t_0)\} + a_{K+1}\{F_{0+}(u) - F_{0+}(t_0)\}\} dG(u)$$

to both sides of (36). This gives

$$\int_{\tau_{nk}}^{s} \{a_k\{\widehat{F}_{nk}(u) - F_{0k}(t_0)\} + a_{K+1}\{\widehat{F}_{n+}(u) - F_{0+}(t_0)\}\} dG(u)$$
$$(37) \qquad \leq \int_{[\tau_{nk},s)} \{a_k\{\delta_k - F_{0k}(t_0)\} + a_{K+1}\{\delta_+ - F_{0+}(t_0)\}\} d\mathbb{P}_n(u,\delta)$$
$$+ R'_{nk}(\tau_{nk}, s),$$

where equality holds if $s$ is a jump point of $\widehat{F}_{nk}$, and where

$$R'_{nk}(s,t) = R_{nk}(s,t) + \rho_{nk}(s,t),$$

with

$$\rho_{nk}(s,t) = \int_{[s,t)} \{a_k\{F_{0k}(t_0) - F_{0k}(u)\}$$
$$+ a_{K+1}\{F_{0+}(t_0) - F_{0+}(u)\}\} d(G_n - G)(u).$$

Note that $\rho_{nk}(\tau_{nk}, s) = o_p(n^{-2/3})$, uniformly in $s \in [t_0 - m_1 n^{-1/3}, t_0 + m_1 n^{-1/3}]$, using (29) in [8], Lemma 4.9 and $t_0 - \tau_{nk} = O_p(n^{-1/3})$ by [8], Corollary 4.19. Hence, the remainder term $R'_{nk}$ in (37) is of the same order as $R_{nk}$. Next, consider (37), and write $\int_{[\tau_{nk},s)} = \int_{[\tau_{nk},t_0)} + \int_{[t_0,s)}$, let $s = t_0 + n^{-1/3}t$, and multiply by $n^{2/3}/g(t_0)$. This yields

$$(38) \quad c_{nk} + a_k \bar{H}_{nk}(t) + a_{K+1} \bar{H}_{n+}(t) \leq R_{nk}^{\mathrm{loc}}(t) + a_k V_{nk}^{\mathrm{loc}}(t-) + a_{K+1} V_{n+}^{\mathrm{loc}}(t-),$$

where equality holds if $t$ is a jump point of $\widehat{F}_{nk}^{\mathrm{loc}}$ and where

$$(39) \quad R_{nk}^{\mathrm{loc}}(t) = \{n^{2/3}/g(t_0)\} R'_{nk}(\tau_{nk}, t_0 + n^{-1/3}t), \qquad k = 1, \ldots, K.$$

Note that $R_{nk}^{\mathrm{loc}}(t) = o_p(1)$ uniformly in $t \in [-m_1, m_1]$, using again that $t_0 - \tau_{nk} = O_p(n^{-1/3})$. Moreover, note that $R_{nk}^{\mathrm{loc}}$ is left-continuous. We now remove the random variables $c_{nk}$ by solving the following system of equations for $H_1, \ldots, H_K$:

$$c_{nk} + a_k \bar{H}_{nk}(t) + a_{K+1} \bar{H}_{n+}(t) = a_k H_{nk}(t) + a_{K+1} H_{n+}(t), \qquad k = 1, \ldots, K.$$

The unique solution is $H_{nk}(t) = \bar{H}_{nk}(t) + (c_{nk}/a_k) + \sum_{k=1}^{K}(c_{nk}/a_k) \equiv \widehat{H}_{nk}^{\mathrm{loc}}(t)$.
□

DEFINITION 3.3. We define $\widehat{U}_n = (R_n^{\mathrm{loc}}, V_n^{\mathrm{loc}}, \widehat{H}_n^{\mathrm{loc}}, \widehat{F}_n^{\mathrm{loc}})$, where $R_n^{\mathrm{loc}} = (R_{n1}^{\mathrm{loc}}, \ldots, R_{nK}^{\mathrm{loc}})$ with $R_{nk}^{\mathrm{loc}}$ defined by (39), and where $V_n^{\mathrm{loc}}$, $\widehat{H}_n^{\mathrm{loc}}$ and $\widehat{F}_n^{\mathrm{loc}}$ are given in Definition 34. We use the notation $\cdot|[-m,m]$ to denote that processes are restricted to $[-m,m]$.

We now define a space for $\widehat{U}_n|[-m,m]$:

DEFINITION 3.4. For any interval $I$, let $D^{-}(I)$ be the collection of "caglad" functions on $I$ (left-continuous with right limits), and let $C(I)$ denote the collection of continuous functions on $I$. For $m \in \mathbb{N}$, we define the space

$$E[-m,m] = (D^{-}[-m,m])^K \times (D[-m,m])^K \times (C[-m,m])^K \times (D[-m,m])^K$$
$$\equiv I \times II \times III \times IV,$$

endowed with the product topology induced by the uniform topology on $I \times II \times III$, and the Skorohod topology on $IV$.

PROOF OF THEOREM 1.8. Analogously to the work of [6], proof of Theorem 6.2, on the estimation of convex densities, we first show that $\widehat{U}_n|[-m,m]$ is tight in $E[-m,m]$ for each $m \in \mathbb{N}$. Since $R_{nk}^{\mathrm{loc}}|[-m,m] = o_p(1)$ by Proposition 3.2, it follows that $R_n^{\mathrm{loc}}$ is tight in $(D^{-}[-m,m])^K$ endowed with the uniform topology. Next, note that the subset of $D[-m,m]$ consisting of absolutely bounded nondecreasing functions is compact in the Skorohod topology. Hence, the local rate of convergence of the MLE (see [8], Theorem 4.17) and the monotonicity of $\widehat{F}_{nk}^{\mathrm{loc}}$, $k = 1,\ldots,K$, yield tightness of $\widehat{F}_n^{\mathrm{loc}}|[-m,m]$ in the space $(D[-m,m])^K$ endowed with the Skorohod topology. Moreover, since the set of absolutely bounded continuous functions with absolutely bounded derivatives is compact in $C[-m,m]$ endowed with the uniform topology, it follows that $\bar{H}_n^{\mathrm{loc}}|[-m,m]$ is tight in $(C[-m,m])^K$ endowed with the uniform topology. Furthermore, $V_n^{\mathrm{loc}}|[-m,m]$ is tight in $(D[-m,m])^K$ endowed with the uniform topology, since $V_n^{\mathrm{loc}}(t) \to_d V(t)$ uniformly on compacta. Finally, $c_{n1}, \ldots, c_{nK}$ are tight since each $c_{nk}$ is the difference of quantities that are tight, using that $t_0 - \tau_{nk} = O_p(n^{-1/3})$ by [8], Corollary 4.19. Hence, also $\widehat{H}_n^{\mathrm{loc}}|[-m,m]$ is tight in $(C[-m,m])^K$ endowed with the uniform topology. Combining everything, it follows that $\widehat{U}_n|[-m,m]$ is tight in $E[-m,m]$ for each $m \in \mathbb{N}$.

It now follows by a diagonal argument that any subsequence $\widehat{U}_{n'}$ of $\widehat{U}_n$ has a further subsequence $\widehat{U}_{n''}$ that converges in distribution to a limit

$$U = (0, V, H, F) \in (C(\mathbb{R}))^K \times (C(\mathbb{R}))^K \times (C(\mathbb{R}))^K \times (D(\mathbb{R}))^K.$$

Using a representation theorem (see, e.g., [2], [15], Representation Theorem 13, page 71, or [17], Theorem 1.10.4, page 59), we can assume that $\widehat{U}_{n''} \to_{\text{a.s.}} U$. Hence, $F = H'$ at continuity points of $F$, since the derivatives of a sequence of convex functions converge together with the convex functions at points where the limit has a continuous derivative. Proposition 3.2 and the continuous mapping theorem imply that the vector $(V, H, F)$ must satisfy

$$\inf_{[-m,m]} \{a_k V_k(t) + a_{K+1} V_+(t) - a_k H_k(t) - a_{K+1} H_+(t)\} \geq 0,$$

$$\int_{-m}^{m} \{a_k V_k(t) + a_{K+1} V_+(t) - a_k H_k(t) - a_{K+1} H_+(t)\} \, dF_k(t) = 0,$$

for all $m \in \mathbb{N}$, where we replaced $V_k(t-)$ by $V_k(t)$, since $V_1, \ldots, V_K$ are continuous.

Letting $m \to \infty$, it follows that $H_1, \ldots, H_K$ satisfy conditions (i) and (ii) of Theorem 1.7. Furthermore, Theorem 1.7(iii) is satisfied since $t_0 - \tau_{nk} = O_p(n^{-1/3})$ by [8], Corollary 4.19. Hence, there exists a $K$-tuple of processes $(H_1, \ldots, H_K)$ that satisfies the conditions of Theorem 1.7. This proves the existence part of Theorem 1.7. Moreover, Proposition 2.16 implies that there is only one such $K$-tuple. Thus, each subsequence converges to the same limit $H = (H_1, \ldots, H_K) = (\widehat{H}_1, \ldots, \widehat{H}_K)$ defined in Theorem 1.8. In particular, this implies that $\widehat{F}_n^{\text{loc}}(t) = n^{1/3}(\widehat{F}_n(t_0 + n^{-1/3}t) - F_0(t_0)) \to_d \widehat{F}(t)$ in the Skorohod topology on $(D(\mathbb{R}))^K$. □

**4. Simulations.** We simulated 1000 data sets of sizes $n = 250$, 2500 and 25,000, from the model given in Example 2.3. For each data set, we computed the MLE and the naive estimator. For computation of the naive estimator, see [1], pages 13–15 and [9], pages 40–41. Various algorithms for the computation of the MLE are proposed by [10, 11, 12]. However, in order to handle large data sets, we use a different approach. We view the problem as a bivariate censored data problem, and use a method based on sequential quadratic programming and the support reduction algorithm of [7]. Details are discussed in [13], Chapter 5. As convergence criterion we used satisfaction of the characterization in [8], Corollary 2.8, within a tolerance of $10^{-10}$. Both estimators were assumed to be piecewise constant, as discussed in the beginning of Section 3.

It was suggested by [12] that the naive estimator can be improved by suitably modifying it when the sum of its components exceeds 1. In order to investigate this idea, we define a "scaled naive estimator" $\widetilde{F}_{nk}^s$ by

$$\widetilde{F}_{nk}^s(t) = \begin{cases} \widetilde{F}_{nk}(t), & \text{if } \widetilde{F}_{n+}(s_0) \leq 1, \\ \widetilde{F}_{nk}(t)/\widetilde{F}_{n+}(s_0), & \text{if } \widetilde{F}_{n+}(s_0) > 1, \end{cases}$$

for $k = 1, \ldots, K$, where we take $s_0 = 3$. Note that $\widetilde{F}_{n+}^s(t) \leq 1$ for $t \leq 3$. We also defined a "truncated naive estimator" $\widetilde{F}_{nk}^t$. If $\widetilde{F}_{n+}(T_{(n)}) \leq 1$, then $\widetilde{F}_{nk}^t \equiv \widetilde{F}_{nk}$ for all $k = 1, \ldots, K$. Otherwise, we let $s_n = \min\{t : \widetilde{F}_{n+}(t) > 1\}$ and define

$$\widetilde{F}_{nk}^t(t) = \begin{cases} \widetilde{F}_{nk}(t), & \text{for } t < s_n, \\ \widetilde{F}_{nk}(t) + \alpha_{nk}, & \text{for } t \geq s_n, \end{cases}$$

where

$$\alpha_{nk} = \frac{\widetilde{F}_{nk}(s_n) - \widetilde{F}_{nk}(s_n-)}{\widetilde{F}_{n+}(s_n) - \widetilde{F}_{n+}(s_n-)} \{1 - \widetilde{F}_{n+}(s_n-)\},$$

for $k = 1, \ldots, K$. Note that $\widetilde{F}_{n+}^t(t) \leq 1$ for all $t \in \mathbb{R}$.

We computed the mean squared error (MSE) of all estimators on a grid with points $0, 0.01, 0.02, \ldots, 3.0$. Subsequently, we computed relative MSEs by dividing the MSE of the MLE by the MSE of each estimator. The results are shown in Figure 3. Note that the MLE tends to have the best MSE, for all sample sizes and for all values of $t$. Only for sample size 250 and small values of $t$, the scaled naive estimator outperforms the other estimators; this anomaly is caused by the fact that this estimator is scaled down so much that it has a very small variance. The difference between the MLE and the naive estimators is most pronounced for large values of $t$. This was also observed by [12], and they explained this by noting that only the MLE is guaranteed to satisfy the constraint $F_+(t) \leq 1$ at large values of $t$. We believe that this constraint is indeed important for small sample sizes, but the theory developed in this paper indicates that it does not play any role asymptotically. Asymptotically, the difference can be explained by the extra term $(a_{K+1}/a_k)\{V_+ - \widehat{H}_+\}$ in the limiting process of the MLE (see Proposition 2.4), since the factor $a_{K+1}/a_k = F_{0k}(t)/F_{0,K+1}(t)$ is increasing in $t$.

Among the naive estimators, the truncated naive estimator behaves better than the naive estimator for sample sizes 250 and 2500, especially for large values of $t$. However, for sample size 25,000 we can barely distinguish the three naive estimators. The latter can be explained by the fact that all versions of the naive estimator are asymptotically equivalent for $t \in [0, 3]$, since consistency of the naive estimator ensures that $\lim_{n \to \infty} \widetilde{F}_{n+}(3) \leq 1$ almost surely. On the other hand, the three naive estimators are clearly less efficient than the MLE for sample size 25,000. These results support our theoretical finding that the form of the likelihood (and not the constrained $F_+ \leq 1$) causes the different asymptotic behavior of the MLE and the naive estimator.

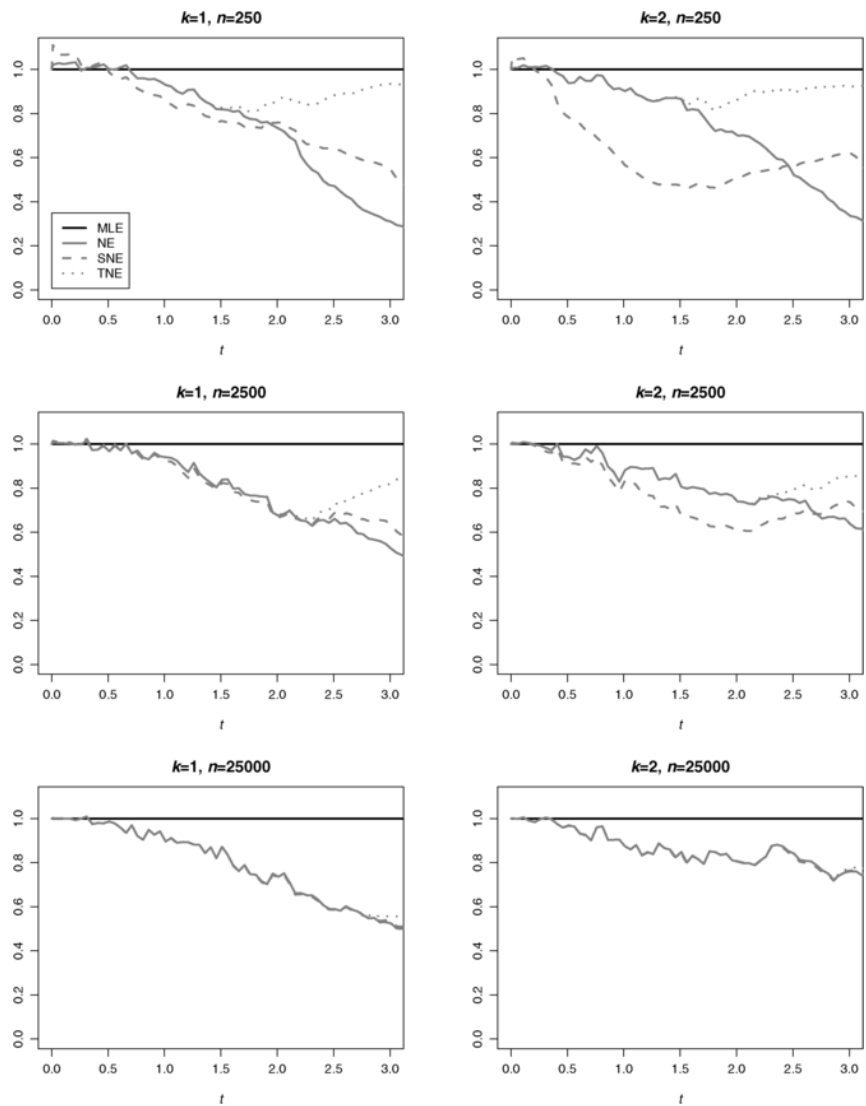

Fig. 3. *Relative MSEs, computed by dividing the MSE of the MLE by the MSE of the other estimators. All MSEs were computed over* 1000 *simulations for each sample size, on the grid* $0, 0.01, 0.02, \ldots, 3.0$.

Finally, we note that our simulations consider estimation of $F_{0k}(t)$, for $t$ on a grid. Alternatively, one can consider estimation of certain smooth functionals of $F_{0k}$. The naive estimator was suggested to be asymptotically efficient for this purpose [12], and [14], Chapter 7, proved that the same is true for the MLE. A simulation study that compares the estimators in this setting is presented in [14], Chapter 8.2.

## 5. Technical proofs.

PROOF LEMMA 2.11. Let $k \in \{1, \ldots, K\}$ and $j \in \mathbb{N} = \{0, 1, \ldots\}$. Note that for $M$ large, we have for all $w \leq j + 1$:

$$C(s_{jM} - w \vee (s_{jM} - w)^{3/2}) \leq \tfrac{1}{2} b(s_{jM} - w)^2.$$

Hence, the probability in the statement of Lemma 2.11 is bounded above by

$$P\bigg\{\sup_{w \leq j+1} \bigg\{\int_w^{s_{jM}} dS_k(u) - \tfrac{1}{2} b(s_{jM} - w)^2\bigg\} \geq 0\bigg\}.$$

In turn, this probability is bounded above by

$$(40) \qquad \sum_{q=0}^{\infty} P\bigg\{\sup_{w \in (j-q, j-q+1]} \int_w^{s_{jM}} dS_k(u) \geq \lambda_{kjq}\bigg\},$$

where $\lambda_{kjq} = b(s_{jM} - (j - q + 1))^2/2 = b(Mv(j) + q - 1)^2/2$.

We write the $q$th term in (40) as

$$P\bigg(\sup_{w \in [j-q, j-q+1)} S_k(s_{jM} - w) \geq \lambda_{kjq}\bigg)$$

$$\leq P\bigg(\sup_{w \in [0, Mv(j)+q)} S_k(w) \geq \lambda_{kjq}\bigg) = P\bigg(\sup_{w \in [0,1)} S_k(w) \geq \frac{\lambda_{kjq}}{\sqrt{Mv(j)+q}}\bigg)$$

$$\leq P\bigg(\sup_{w \in [0,1]} B_k(w) \geq \frac{\lambda_{kjq}}{b_k\sqrt{Mv(j)+q}}\bigg)$$

$$\leq 2P\bigg(N(0,1) \geq \frac{\lambda_{kjq}}{b_k\sqrt{Mv(j)+q}}\bigg) \leq 2b_{kjq} \exp\bigg(-\frac{1}{2}\bigg(\frac{\lambda_{kjq}}{b_k\sqrt{Mv(j)+q}}\bigg)^2\bigg),$$

where $b_k$ is the standard deviation of $S_k(1)$ and $b_{kjq} = b_k\sqrt{Mv(j)+q}/(\lambda_{kjq} \times \sqrt{2\pi})$, and $B_k(\cdot)$ is standard Brownian motion. Here we used standard properties of Brownian motion. The second to last inequality is given in, for example, [16], (6), page 33, and the last inequality follows from Mills' ratio (see [3], (10)). Note that $b_{kjq} \leq d$ all $j \in \mathbb{N}$, for some $d > 0$ and all $M > 3$. It follows that (40) is bounded above by

$$\sum_{q=0}^{\infty} d \exp\bigg(-\frac{1}{2}\bigg(\frac{\lambda_{kjq}}{b_k\sqrt{Mv(j)+q}}\bigg)^2\bigg) \approx \sum_{q=0}^{\infty} d \exp\bigg(-\frac{1}{2}\frac{(Mv(j)+q)^3}{b_k^2}\bigg),$$

which in turn is bounded above by $d_1 \exp(-d_2(Mv(j))^3)$, for some constants $d_1$ and $d_2$, using $(a+b)^3 \geq a^3 + b^3$ for $a, b \geq 0$. $\square$

PROOF OF LEMMA 2.12. This proof is completely analogous to the proof of Lemma 4.14 of [8], upon replacing $\widehat{F}_{nk}(u)$ by $\widehat{F}_k(u)$, $F_{0k}(u)$ by $\bar{F}_{0k}(u)$, $dG(u)$ by $du$, $S_{nk}(\cdot)$ by $S_k(\cdot)$, $\tau_{nkj}$ by $\tau_{kj}$, $s_{njM}$ by $s_{jM}$, and $A_{njM}$ by $A_{jM}$. The only difference is that the second term on the right-hand side of equation (69) in [8], vanishes, since this term comes from the remainder term $R_{nk}(s,t)$, and we do not have such a remainder term in the limiting characterization given in Proposition 3.2. $\square$

P. GROENEBOOM
DEPARTMENT OF MATHEMATICS
DELFT UNIVERSITY OF TECHNOLOGY
MEKELWEG 4
2628 CD DELFT
THE NETHERLANDS
E-MAIL: p.groeneboom@ewi.tudelft.nl

M. H. MAATHUIS
SEMINAR FÜR STATISTIK
ETH ZÜRICH
CH-8092 ZÜRICH
SWITZERLAND
E-MAIL: maathuis@stat.math.ethz.ch

J. A. WELLNER
DEPARTMENT OF STATISTICS
UNIVERSITY OF WASHINGTON
SEATTLE, WASHINGTON 98195
USA
E-MAIL: jaw@stat.washington.edu